\documentclass[11pt]{article}
\usepackage{amsfonts,latexsym,amsmath,amscd,geometry}
\usepackage{hyperref}
\usepackage{mathtools}
\usepackage{amssymb}
\usepackage[normalem]{ulem}
\usepackage{latexsym}
\usepackage{xcolor}

\makeatletter \@addtoreset{equation}{section} \makeatother

\newcommand \nc{\newcommand}
\newtheorem{theorem}{Theorem}[section]
\newtheorem{lemma}[theorem]{Lemma}

\newtheorem{definition}[theorem]{Definition}

\nc{\ba}{\begin{array}}\nc{\ea}{\end{array}}
\nc{\be}{\begin{eqnarray}}\nc{\ee}{\end{eqnarray}}
\nc{\beq}{\begin{equation}}\nc{\eeq}{\end{equation}}
\nc{\bex}{\begin{eqnarray*}}\nc{\eex}{\end{eqnarray*}}
\nc{\btm}{\begin{theorem}} \nc{\etm}{\end{theorem}}
\nc{\blm}{\begin{lemma}} \nc{\elm}{\end{lemma}}
\let\isout\sout
\renewcommand{\sout}[1]{\ifmmode\text{\isout{\ensuremath{#1}}}\else\isout{#1}\fi}
\nc{\R}{\mathbb{R}}  
\def\pf{\noindent{\bf Proof.\quad}}
\def\endpf{\hfill$\Box$}
\def\les{\lesssim}

\def\x{{\bf x}}
\def\y{{\bf y}}
\def\vv{{v}}
\def\hv{\widehat{{v}}}
\def\ve{\varepsilon}

\def\pa{\partial}

\def\f{{\bf f}}
\def\N{\mathcal{N}}
\DeclareMathOperator{\loc}{loc}
\DeclareMathOperator{\dist}{dist}
\DeclareMathOperator{\dv}{div}
\begin{document}

\title{Weak compactness of simplified nematic liquid flows in 2D}

\author{
Hengrong Du\footnote{Department of Mathematics, Purdue University, West Lafayette, IN 47906, USA} 
\quad Tao Huang\footnote{Department of Mathematics, Wayne State University, Detroit, MI 48202, USA}  \quad Changyou Wang\footnote{Department of Mathematics, Purdue University, West Lafayette, IN 47906, USA}\\
}

\maketitle

\begin{abstract}
For any bounded, smooth domain $\Omega\subset \R^2$, 
we will establish the weak compactness property of solutions to the simplified Ericksen-Leslie system for both uniaxial and biaxial nematics, and the convergence of weak solutions of  the Ginzburg-Landau type nematic liquid crystal flow to a weak solution of  the simplified Ericksen-Leslie system as the parameter tends to zero. This is based on the compensated compactness property of the Ericksen stress tensors, which is obtained by the $L^p$-estimate of the Hopf differential for the Ericksen-Leslie system and the Pohozaev type argument for the Ginzburg-Landau type 
nematic liquid crystal flow.

\end{abstract}

\section{Introduction}

Let $\Omega\subset\R^2$ be a bounded domain with smooth boundary, and $\N\subset \R^L$ (for $L\geq 2$) be a smooth compact Riemannian manifold without boundary, and $0<T\le\infty$. 
 We formulate a generalized form of simplified Ericksen-Leslie system of nematic liquid crystals in which the director
 field takes values in $\N$: 
\begin{equation}
\begin{cases}
 u_t+u\cdot\nabla u-\Delta u+\nabla P=-\nabla\cdot\big(\nabla v\odot\nabla v\big),\\
\nabla \cdot u=0, \\
 v_t+u\cdot\nabla v=\Delta v+A(v)(\nabla v, \nabla v),
\end{cases}
\ {\rm{in}}\ \Omega\times (0, T),
\label{NEL}
\end{equation}
where $(u(\x,t), v(\x,t), P(\x,t)):\Omega\times (0,T)\rightarrow \R^2\times\N\times \R$ represents the fluid velocity field, the orientation director field of nematic material (into a general Riemannian manifold), and the pressure function respectively, $\big(\nabla v\odot\nabla v\big)_{ij}=\nabla_{\x_i}v\cdot\nabla_{\x_j} v$ for $i,j=1, 2$ represents the Ericksen-Leslie stress tensor, and $A(y)(\cdot, \cdot)$ is the second fundamental form of $\N$ at the point
$y\in\N$.  

The generalized system \eqref{NEL} covers the two important cases in nematic liquid crystals:
\begin{enumerate}
\item[(1)] For $\N=\mathbb S^2$, the system \eqref{NEL} becomes the simplified, uniaxial Ericksen-Leslie system first proposed by \cite{lin1989nonlinear} 
\begin{equation}\label{ELlc}
\begin{cases}
\pa_t u+u\cdot\nabla u-\Delta u+\nabla P=-\nabla\cdot\big(\nabla d\odot\nabla d\big),\\
\nabla \cdot u=0, \\
\pa_t d+u\cdot\nabla d=\Delta d+|\nabla d|^2d,
\end{cases}
\end{equation} 
 for $(u(\x,t), d(\x,t), P(\x,t)):\Omega\times (0,T)\rightarrow \R^2\times\mathbb S^2\times \R$. 
In dimension two, the existence of a unique global weak solution has been proved in \cite{lin2010liquid} \cite{lin2010uniqueness}, which satisfies the energy inequality and has at most finitely many singular times,
see also \cite{Hong2011}. Very recently, the authors in \cite{LLWWZ19} have constructed example of singularity at finite time. In dimension three, a global weak solution has been constructed in \cite{linwang2016weakthreedim} with initial data $d_0\in \mathbb S^2_+$. Examples of finite time singularity have been constructed by \cite{Huang2016FiniteTS}. Interested readers
can consult the survey article \cite{lin2014recent} and the references therein.

\item[(2)] 
For
$$
\mathcal N=\big\{({\bf y}_1, {\bf y}_2)\in \mathbb S^2\times \mathbb S^2~\big|~{\bf y}_1\cdot{\bf y}_2=0\big\}
\subset \mathbb R^6,
$$
let $v(\x, t)=(n(\x, t), m(\x, t)):\Omega\times (0,T)\rightarrow\mathbb S^2\times\mathbb S^2$ with $n\cdot m=0$. Then the system \eqref{NEL} becomes the biaxial, Ericksen-Leslie system
\begin{equation}\label{blc}
\begin{cases}
\pa_t u+u\cdot\nabla u-\Delta u+\nabla P=-\nabla\cdot\big(\nabla n\odot\nabla n+\nabla m\odot\nabla m\big),\\
\nabla \cdot u=0, \\
\pa_t n+u\cdot\nabla n=\Delta n+|\nabla n|^2n+\langle\nabla n,\nabla m\rangle m\\
\pa_t m+u\cdot\nabla m=\Delta m+|\nabla m|^2m+\langle\nabla m,\nabla n\rangle n\\
n\cdot m=0
\end{cases}
\ {\rm{in}}\ \Omega\times (0,T).
\end{equation}
This is a simplified version of the hydrodynamics of biaxial nematics model proposed by Grovers and
Vertogen \cite{GoversVertogen1984a, GoversVertogen1984b, GoversVertogen1985}. 
In  dimensional two, the existence of a unique global weak solution has recently been shown in \cite{LinLiWang20},
which is smooth off at most finitely many singular times. 

\end{enumerate}

A strategy to construct a weak solution of \eqref{NEL} and \eqref{inbdy} is to consider a Ginzburg-Landau approximated system (cf. \cite{linliu1995nonparabolic},\cite{linliu1996partialregularity}). 
More precisely,  for any $\delta>0$ set the $\delta$-neighborhood of $\N$ by
$$
\mathcal N_\delta=\Big\{{\bf y}\in \mathbb R^L~\big|~ \dist({\bf y}, \mathcal N)<\delta\Big\},
$$
where $\dist (\y, \mathcal N)$ is the distance from $\y$ to $\mathcal N$. Let $\Pi_{\mathcal N}:\mathcal N_\delta\rightarrow \mathcal N$ be the nearest point projection map. There exists $\delta_{\N}=\delta(\mathcal N)>0$ such that $\dist (\y, \mathcal N)$ and 
$\Pi_{\N}$ are smooth in $\N_{2\delta_{\N}}$. Let $\chi(s)\in C^\infty([0,\infty))$ be a monotone increasing function such that 
$$
\chi(s)=
\left\{
\begin{array}{ll}
s, \quad & \mbox{if } 0\leq s
\leq \delta_{\N}^2,\\
4 \delta_{\N}^2, \quad &\mbox{if }  s
\geq 4\delta_{\N}^2.
\end{array}
\right.
$$
Consider the following Ginzburg-Landau energy functional for the director $v$
$$
E_{\varepsilon}(u, {v})
=\int_{\Omega}\big(\frac12|\nabla {v}|^2+\frac{1}{\varepsilon^2}\chi\big(\dist^2({v}, \mathcal N)\big)\big).
$$
Then the corresponding Ginzburg-Landau approximated system of \eqref{NEL} can be written as
\begin{equation}\label{NGL}
\begin{cases}
u_t+u\cdot\nabla u-\Delta u+\nabla P=-\nabla\cdot\big(\nabla {v}\odot\nabla {v}\big),\\
\nabla \cdot u=0, \\
\displaystyle{v}_t+u\cdot\nabla {v}=\Delta {v}-\frac{1}{\ve^2}\chi'\big(\dist^2({v}, \mathcal N)\big)\frac{d}{d{v}}\big(\dist^2({v}, \mathcal N)\big).
\end{cases}
\end{equation}

The main purpose of this paper is to study the weak compactness of  solutions to the simplified Ericksen-Leslie system \eqref{NEL} and convergence of solutions of the Ginzburg-Landau  approximation \eqref{NGL} to the
simplified Ericksen-Leslie system \eqref{NEL}.  For this purpose, we will consider 
the following initial and boundary condition 
\beq\label{inbdy}
(u, v)\,|\,_{\partial_p Q_T}=(u_0, v_0)
\eeq
where $Q_T=\Omega\times(0,T)$ and $\partial_p Q_T=\big(\Omega\times\{t=0\}\big)\cup\big(\partial \Omega\times[0,T]\big)$ is the parabolic boundary of $Q_T$. 
We assume that 
\beq\label{uoc}
u_0\,\big|\,_{\partial\Omega}=0,\quad v_0(x)\in \N\ \ \mbox{for  a.e. }\ x\in \Omega,
\eeq
and  introduce the notations
\begin{equation*}
  \begin{split}
  {\bf H}=\text{closure of }C_0^\infty\left(\Omega, \R^2\right)\cap \left\{ f~\big|~\nabla\cdot f=0\right\} \text{ in }L^2\left(\Omega,\R^2\right),\\
  {\bf J}=\text{closure of }C_0^\infty\left(\Omega, \R^2\right)\cap \left\{ f~\big|~\nabla\cdot f=0 \right\}\text{ in }H_0^1\left(\Omega, \R^2\right), \\
  H^1(\Omega, \N)=\left\{ f\in H^1(\Omega, \R^L)~\big|~ f(x)\in \N \text{ a.e. }x\in \Omega \right\}.
  \end{split}
\end{equation*}
We also assume that 
\begin{equation}  \label{initialRegularity}
  u_0\in {\bf H},\ \ v_0\in H^1(\Omega, \N). 
\end{equation}
Recall the definition of weak solutions of \eqref{NEL}. 

\begin{definition}\label{defwsol}
A pair of maps $u\in L^{\infty}([0,T], {\bf H})\cap L^2([0,T], {\bf J})$ and $v\in L^2([0,T], H^1(\Omega, \N))$ is called a weak solution to initial and boundary problem \eqref{NEL}, \eqref{inbdy}-\eqref{initialRegularity}, if 
  \begin{equation}
  \begin{split}
    &-\int_{Q_T}\left\langle u, \xi' \varphi  \right\rangle+\int_{Q_T} \left\langle u\cdot \nabla u, \xi\varphi \right\rangle+\left\langle \nabla u, \xi \nabla\varphi \right\rangle\\
    &\ =-\xi(0)\int_\Omega\left\langle u_0, \varphi \right\rangle+\int_{Q_T} \left\langle \nabla v\odot\nabla v, \xi \nabla \varphi \right\rangle, \\
    &-\int_{Q_T} \left\langle  v, \xi'\phi\right\rangle+\int_{Q_T}\left\langle u\nabla v, \xi\phi \right\rangle
    +\left\langle \nabla v, \xi \nabla\phi \right\rangle\\
    &\ =-\xi(0)\int_\Omega\left\langle v_0, \phi \right\rangle
    +\int_{Q_T} \left\langle A(v)(\nabla v, \nabla v), \xi\phi \right\rangle,
    \label{weakdef}
    \end{split}
  \end{equation}
  for any $\xi\in C^\infty([0,T])$ with $\xi(T)=0$, $\varphi\in {\bf J}$ and $\phi \in H_0^1(\Omega, \R^3)$. Moreover, $(u, v)|_{\partial \Omega}=(u_0, v_0)$ in the sense of trace. The notion of a weak solution to 
the system \eqref{NGL} can be defined similarly.    
\end{definition}

Our first main theorem concerns the convergence of weak solutions of the system \eqref{NGL}
to the system \eqref{NEL} as $\varepsilon\rightarrow 0$. We remark that the existence of weak solutions
to \eqref{NGL} has been established by \cite{linliu1995nonparabolic, linliu1996partialregularity} for $\N=\mathbb S^2$ by the Galerkin method, which can be easily adapted to handle the case that $\N$ is a compact Riemannian manifold.

\begin{theorem}  \label{thm3}
  For $\ve>0$, let $(u^\ve, {v}^\ve)$ be a sequence of weak solutions to the Ginzburg-Landau approximated system \eqref{NGL} with the initial and boundary condition \eqref{inbdy}-\eqref{initialRegularity}. Then there exists a weak solution $(u, v)$ of \eqref{NEL} with the initial and boundary condition \eqref{inbdy}-\eqref{initialRegularity} such that, after passing to subsequences, 
\begin{equation}\notag
 u^\ve\rightharpoonup u\ \text{ in }\ L^2 ([0,T], H^1(\Omega)) , \quad 
      {v}^\ve \rightharpoonup v\ \text{ in }\ L^2 ([0,T], H^1(\Omega)) .
  \end{equation}
In particular, the initial and boundary problem \eqref{NEL} and \eqref{inbdy}-\eqref{initialRegularity} admits at least one weak solution $u\in L^{\infty}([0,T], {\bf H})\cap L^2([0,T], {\bf J})$ and $v\in L^2([0,T], H^1(\Omega, \N))$.
\end{theorem}

We would like to mention that when $\N=\mathbb S^2$, the convergence of solutions of system \eqref{NGL} to the system \eqref{ELlc} has recently been proved in two dimensional torus $T^2$ by Kortum in an interesting article \cite{kortum2019concentrationcancellation}. In order to deal with convergence of the most difficult terms $\nabla d_\varepsilon\odot\nabla d_\varepsilon$ in the limit process, Kortum employed the concentration-cancellation method for the Euler equation developed by DiPerna and Majda \cite{DipernaMajda1988concentration} (see also \cite{MajdaBertozzi2002Vorticity}).  Thanks to the rotational covariance of $\nabla d_\varepsilon\odot\nabla d_\varepsilon$, the test functions can be taken to a function of periodic one spatial variable 
ensuring the weak convergence of $\nabla d_\varepsilon\odot\nabla d_\varepsilon$ to $\nabla d\odot\nabla d$.

In this paper, we make some new observations on the Ericksen stress tensor $\nabla v\odot \nabla v$, which is flexible enough to handle any smooth domain $\Omega\subset \R^2$.  Namely, by adding $-\frac{1}{2}|\nabla v^\ve|^2\mathbb I_2$ to $\nabla v^\ve\odot \nabla v^\ve$, where $\mathbb I_2$ is the $2\times2$ identity matrix, 
we have
$$
\nabla v^\ve\odot\nabla v^\ve-\frac{1}{2}|\nabla v^\ve|^2\mathbb I_2
=\frac12\left(
\begin{array}{ll}
|\partial_x v^\ve|^2-|\partial_y v^\ve|^2, &2\langle\partial_x v^\ve,\partial_y v^\ve\rangle\\
2\langle\partial_x v^\ve,\partial_y v^\ve\rangle, &|\partial_y v^\ve|^2-|\partial_x v^\ve|^2
\end{array}
\right).
$$
This is a matrix whose components constitute the Hopf differential of map $v^\ve$, which are 
$|\partial_x v^\ve|^2-|\partial_y v^\ve|^2$ and $\langle\partial_x v^\ve,\partial_y v^\ve\rangle$. 
Since $v^\epsilon$ is either an approximated harmonic map to $\N$ or a Ginzburg-Landau
approximated harmonic map, we can develop its compensated compactness property
by the Pohozaev type argument.

\medskip
As a byproduct of the proof of Theorem \ref{thm3}, we obtain the following compactness for a sequence of weak solutions to the system \eqref{NEL}.

\begin{theorem}\label{thm1}
 Let $(u^k, v^k):\Omega\times (0,T)\rightarrow \R^2\times\N$ be a sequence of weak solutions to \eqref{NEL}, 
 along with the initial and boundary condition $(u^k_0,v^k_0)$ satisfying \eqref{uoc}, such that
 \beq\label{kengbd}
\begin{split}
\sup\limits_{k\geq 1} \left\{\int_{Q_t}\big(|u^k|^2+|\nabla v^k|^2\big)+\int_{Q_t}\big(|\nabla u^k|^2+|v^k_t+u^k\cdot\nabla v^k|^2\big)\right\}< \infty,
\end{split}
\eeq
Furthermore, if we assume that 
\beq\notag
(u^k_0,v^k_0)\rightharpoonup (u_0, v_0)\ \ \mbox{in }\ \ L^2(\Omega)\times H^1(\Omega),
\eeq
then there exists a weak solution $(u, v)$ of \eqref{NEL} with the initial and boundary condition $(u_0, v_0)$
such that, after passing to subsequences, 
\begin{equation}\notag
  u^k\rightharpoonup u \text{ in } L^2 ([0,T], H^1(\Omega)),  \qquad 
      {v}^k \rightharpoonup v \text{ in } L^2 ([0,T], H^1(\Omega)).
\end{equation}
\end{theorem}

Since the system \eqref{NEL} possesses the geometric structure, i.e., 
$$
A(v^k)(\cdot, \cdot) \perp  T_{v^k}\N
$$
where $T_{v^k}\N$ is the tangent space of $\N$ at $v^k$, 
we can show the weak convergence of $|\partial_x v^k|^2-|\partial_y v^k|^2$ and $\langle\partial_x v^k,\partial_y v^k\rangle$ by utilizing the $L^p$-estimate, $1<p<2$, of the Hopf differential of $v^k$.

\section{Estimates on inhomogeneous Ginzburg-Landau equations}

In this section, we will consider the inhomogeneous Ginzburg-Landau equation
\beq\label{appGLb}
\Delta {v}^\ve-\frac{1}{\ve^2}\chi'\big(\dist^2({v}^\ve, \mathcal N)\big)\frac{d}{d{v}}\big(\dist^2({v}^\ve, \mathcal N)\big)=\tau_\ve
\quad \mbox{in }\ \ \Omega.
\eeq
Suppose 
\begin{equation}
  \sup_{0<\ve\le 1} \mathcal{E}_\ve({v}^\ve)=\int_\Omega\left( \frac{1}{2}|\nabla {v}^\ve|^2+\frac{1}{\ve^2}\chi(\dist^2(\vv^\ve, \mathcal{N})) \right)\le \Lambda_1<\infty, 
  \label{eqn:Lambda1}
\end{equation}
and
\begin{equation}
  \sup_{0<\ve\le 1}\left\|\tau^\ve\right\|_{L^2(\Omega)}\le \Lambda_2<\infty.
  \label{eqn:Lambda2}
\end{equation}
Assume that there exist ${v}\in H^1(\Omega, \mathcal{N})$ and $\tau\in L^2(\Omega,\mathbb R^L)$ 
such that 
\begin{equation}\notag
    \tau^\ve\rightarrow \tau \text{ in }L^2(\Omega), \quad
    \vv^\ve\rightharpoonup \vv  \text{ in }H^1(\Omega).
\end{equation}
Then we have
\begin{lemma}\label{lemma:biGL}
  There exists $\delta_0>0$  such that if $\vv^\ve\in H^1(\Omega, \R^L)$ is a family of solutions to \eqref{appGLb} satisfying \eqref{eqn:Lambda1} and \eqref{eqn:Lambda2}, and for $\x_0\in\Omega$
  and $0<r_0<\dist(x_0,\partial\Omega)$,
\beq\label{engdelta}
\sup_{0<\ve\le 1}\int_{B_{r_0}(\x_0)}\big(\frac12|\nabla {v}^\ve|^2+\frac{1}{\varepsilon^2}\chi\big(\dist^2({v}^\ve, \mathcal N)\big)\big)\leq \delta_0^2,
\eeq  
then there exists an approximated harmonic map ${v}\in H^1(B_{\frac{r_0}{4}}(\x_0),\mathcal N)$ with tension filed $\tau$, i.e, 
\beq
\Delta {v}+A({v})(\nabla {v},\nabla {v})=\tau,
\eeq
such that as $\epsilon\rightarrow 0$, 
\beq
{v}^\ve\rightarrow {v}\quad \mbox{in }\ H^1(B_{\frac{r_0}4}(\x_0)),\ \ {\rm{and}}\ 
\ \frac{1}{\varepsilon^2}\chi\big(\dist^2({v}^\ve, \mathcal N)\big)\rightarrow 0\ \ {\rm{in}}\ \
L^1(B_{\frac{r_0}4}(\x_0)).
\eeq
\end{lemma}

\pf For any fixed $\x_1 \in B_{\frac{r_0}2}(\x_0)$ 
and $0<\ve\le \frac{r_0}{2}$, define $\hv^\ve(\x)=\vv^\ve(\x_1+\ve \x):B_1(0)\to \R^L$. Then we have
\begin{equation*}
  \Delta \hv^\ve=\chi'(\dist^2(\hv^\ve, \mathcal{N}))\frac{d}{d\vv}(\dist^2(\hv^\ve, \mathcal{N}))+\widehat{\tau}^\ve \text{ in }B_1(0),
\end{equation*}
where $\widehat{\tau}^\ve(\x)=\ve^2 \tau^\ve(\x_1+\ve \x)$. Since 
\begin{align*}
  \left\|\Delta \hv^\ve\right\|_{L^2(B_1(0))}&\le \big\| \chi'(\dist(\hv^\ve, \mathcal{N}))\frac{d}{d\vv}(\dist^2(\hv^\ve, \mathcal{N})) \big\|_{L^2(B_1(0))}+\left\|\widehat{\tau}^\ve\right\|_{L^2(B_1(0))}\\
  &\le C\Big(\int_{\Omega\cap \{\dist(v^\epsilon, \N)\le 2\delta_{\N}\}}
  |\dist(v^\epsilon, \N)|^2\Big)^\frac12+\ve\left\|\tau^\ve\right\|_{L^2(\Omega)}\le C+\Lambda_2.
\end{align*}
Thus $\hv^{\ve}\in H^2(B_\frac12)$ and $\left\|\hv^\ve\right\|_{H^2(B_\frac12)}\le C(1+\Lambda_2)$. By Morrey's inequality, we conclude that $\hv^\ve\in C^{\frac12}(B_\frac12)$ and 
\begin{equation*}
  [\hv^{\ve}]_{C^{\frac12}(B_\frac12)}\le C\left\|\hv^\ve\right\|_{H^2(B_\frac12)}\le C(1+\Lambda_2).
\end{equation*}
By rescaling, we get 
\begin{equation*}
  |\hv^\ve(\x)-\hv^\ve(\y)|\le C(1+\Lambda_2)\Big( \frac{|\x-\y|}{\ve} \Big)^{\frac{1}{2}}, \qquad \forall \x, \y\in B_\ve(\x_1).
\end{equation*}
We claim that $\dist(\vv^\ve, \mathcal{N})\le \delta_{\mathcal{N}}$ on $B_{\frac{r_0}2}(\x_0)$. Suppose it were false. Then there exists $\x_1\in B_{\frac{r_0}2}(\x_0)$ such that $\dist(\vv^\ve(\x_1), \mathcal{N})>\delta_{\mathcal{N}}$. Then for any $\theta_0\in (0, 1)$ and $\x\in B_{\theta_0 \ve}(\x_1)$, it holds
\begin{equation*}
  |\vv^\ve(\x)-\vv^\ve(\x_1)|\le C\Big( \frac{|\x-\x_1|}{\ve} \Big)^{\frac12}\le C \theta_0^{\frac12}\le \frac{1}{2}\delta_{\mathcal{N}},
\end{equation*}
provided $\theta_0\le \frac{\delta_{\mathcal{N}}^2}{4C^2}$. It follows that 
\begin{equation*}
  \dist(\vv^\ve(\x), \mathcal{N})\ge \frac{1}{2}\delta_{\mathcal{N}}, \ \forall\x\in B_{\theta_0\ve}(\x_1),
\end{equation*}
so that 
\begin{equation*}
  \int_{B_{\theta_0 \ve}(\x_1)}\frac{1}{\ve^2}\chi(\dist^2(\vv^\ve, \mathcal{N}))\ge \pi \delta_{\mathcal{N}}^2 \theta_0^2. 
\end{equation*}
which contradicts to the assumption that 
\begin{equation*}
  \int_{B_{\theta_0 \ve}(\x_1)}\frac{1}{\ve^2}\chi(\dist^2(\vv^\ve, \mathcal{N}))\le \int_{B_{r_1}(0)}\Big( \frac{1}{2}|\nabla\vv^\ve|^2+\frac{1}{\ve^2}\chi\left(\dist^2(\vv^\ve, \mathcal{N})\right) \Big)\le \delta_0^2
\end{equation*}
for a sufficiently small $\delta_0>0$.

From $\dist(\vv^\ve, \mathcal{N})\le \delta_{\mathcal{N}}$ in $B_{\frac{r_0}2}(\x_0)$, 
we may decompose ${v}^\ve $ into
$$
{v}^\ve=\Pi_{\mathcal N}({v}^\ve) 
+\dist({v}^\ve, \mathcal N)\nu\big(\Pi_{\mathcal N}({v}^\ve)\big):=\omega_\ve+\zeta_\ve\nu_\ve,
$$
so that the equation of ${v_\ve}$ becomes
\beq\label{eqnwdnu}
\Delta \omega_\ve+\Delta \zeta_\ve\nu_\ve+ 2\nabla \zeta_\ve\nabla\nu_\ve+\zeta_\ve\Delta\nu_\ve-\frac{1}{\ve^2}\chi'(\zeta_\ve^2)\nabla_{{v_\ve}}\zeta_\ve^2=\tau_\ve.
\eeq
Multiplying \eqref{eqnwdnu} by $\nu_\ve$, we get
\beq\label{eqnd}
\Delta \zeta_\ve=\langle \nabla \omega_\ve, \nabla\nu_\ve\rangle+\zeta_\ve|\nabla \nu_\ve|^2+\frac{1}{\ve^2}\chi'(\zeta_\ve^2)\langle\nabla_{{v_\ve}}\zeta_\ve^2,\nu_\ve\rangle+\tau^\perp_\ve,
\eeq
where $\tau^\perp_\ve=\langle\tau_\ve, \nu_\ve\rangle$. Plugging $\Delta \zeta_\ve$ into \eqref{eqnwdnu}, we obtain
\beq\label{eqnwnu}
\Delta \omega_\ve+\langle \nabla \omega_\ve,\nabla \nu_\ve\rangle\nu_\ve+\zeta_\ve\big(\Delta \nu_\ve+|\nabla \nu_\ve|^2\nu_\ve\big)+2\langle\nabla\nu_\ve,\nabla \zeta_\ve\rangle=\tau''_\ve,
\eeq
where $\tau''_\ve=\tau_\ve-\tau^\perp_\ve\nu_\ve$. Here we have used the fact
$$
\langle\nabla_{{v_\ve}}\zeta_\ve^2,\nu_\ve\rangle\nu_\ve
=\nabla_{{v_\ve}}\zeta_\ve^2.
$$

Let $\eta\in C^{\infty}_0(B_{\frac{r_0}2}(\x_0), \R)$ be a standard cutoff function of $B_{\frac{3r_0}8}(\x_0)$. 
Since $\dist(\vv^\ve, \mathcal{N})\le \delta_{\mathcal{N}}$, we have that $\chi'(\zeta_\ve^2)=1$ and hence
\beq\label{eqndeta}
\begin{split}
\left(-\Delta+\frac{2}{\ve^2}\right) (\zeta_\ve\eta^2)=&-\zeta_\ve\Delta(\eta^2)-2\nabla \zeta_\ve\nabla(\eta^2)+\langle \nabla \omega_\ve, \nabla(\nu_\ve\eta^2)\rangle-\langle \nabla \omega_\ve, \nu_\ve\nabla(\eta^2)\rangle\\
&+\zeta_\ve\left(|\nabla (\nu_\ve\eta^2)|^2-|\nu_\ve\nabla(\eta^2)|^2\right)+\tau^\perp_\ve\eta^2.
\end{split}
\eeq
Applying the $W^{2,\frac43}$-estimate for $(-\Delta+\frac2{\ve^{2}})$ (see \cite{Krylov2008Sobolev}), we obtain
\beq
\begin{split}
&\|\nabla^2(\zeta_\ve\eta^2)\|_{L^{\frac43}}\\
\les&\|\zeta_\ve\Delta (\eta^2)\|_{L^{\frac43}}+\|\nabla  \zeta_\ve\nabla(\eta^2)\|_{L^{\frac43}}
+\|\nabla \omega_\ve\|_{L^{2}}\|\nabla(\nu_\ve\eta^2)\|_{L^{4}}\\
&+\|\nabla \omega_\ve\|_{L^{2}}+\|\zeta_\ve\|_{L^{\infty}}\|\nabla \nu_\ve\|_{L^{2}}\|\nabla(\nu_\ve\eta^2)\|_{L^{4}}+\|\tau^\perp_\ve\|_{L^{\frac43}}\\
\les& \|\zeta_\ve\|_{L^{\infty}}+\|\nabla \zeta_\ve\|_{L^{2}}+\|\nabla \omega_\ve\|_{L^{2}}(\|\nabla(\nu_\ve\eta^2)\|_{L^{4}}+1)\\
&+\|\zeta_\ve\|_{L^{\infty}}\|\nabla\nu_\ve\|_{L^{2}}\|\nabla(\nu_\ve\eta^2)\|_{L^{4}}+\|\tau_\ve\|_{L^{2}},
\end{split}
\eeq
where $A\les B$ stands for $A\leq CB$ for some universal positive constant $C$.

For $\omega_\ve$, by a similar calculation we obtain 
\beq\label{eqnweta}
\begin{split}
\Delta (\omega_\ve\eta^2)=&-\langle \nabla \omega_\ve,\nabla (\nu_\ve\eta^2)\rangle\nu_\ve+\langle \nabla \omega_\ve, \nu_\ve\nabla(\eta^2)\rangle\nu_\ve\\
&-\zeta_\ve\big[\Delta (\nu_\ve\eta^2)-\nu_\ve\Delta(\eta^2)-2\nabla \nu_\ve\nabla(\eta^2)\big]+\zeta_\ve\big[|\nabla \nu_\ve\eta^2|^2-|\nu_\ve\nabla(\eta^2)|^2\big]\nu_\ve\\
&-2\big[\langle\nabla(\nu_\ve\eta^2),\nabla \zeta_\ve\rangle-\langle\nabla(\eta^2),\nabla \zeta_\ve\rangle\nu_\ve\big]+\tau''_\ve\eta^2+\omega_\ve\Delta(\eta^2)+2\nabla \omega_\ve\nabla (\eta^2).
\end{split}
\eeq
Applying the $W^{2,\frac43}$-estimate, we obtain
\beq
\begin{split}
&\|\nabla^2(\omega_\ve\eta^2)\|_{L^{\frac43}}\\
\les& \|\nabla \omega_\ve\|_{L^{2}}\|\nabla (\nu_\ve\eta^2)\|_{L^{4}}+\|\nabla \omega_\ve\|_{L^{\frac43}}+\|\zeta_\ve\|_{L^{\infty}}\|\Delta(\nu_\ve\eta^2)\|_{L^{\frac43}}\\
&+\|\zeta_\ve\|_{L^{\infty}}\left(1+\|\nabla\nu_\ve\|_{L^{\frac43}}\right)+\|\zeta_\ve\|_{L^{\infty}}\|\nabla(\nu_\ve\eta^2)\|_{L^{2}}\|\nabla(\nu_\ve\eta^2)\|_{L^{4}}\\
&+\|\nabla \zeta_\ve\|_{L^{2}}\|\nabla(\nu_\ve\eta^2)\|_{L^{4}}+\|\nabla \zeta_\ve\|_{L^{\frac43}}+\|\tau_\ve\|_{L^{2}}.
\end{split}
\eeq
Therefore, we conclude that 
\beq
\begin{split}
&\|\nabla^2(\zeta_\ve\eta^2)\|_{L^{\frac43}}+\|\nabla^2(\omega_\ve\eta^2)\|_{L^{\frac43}}\\
\les&\|\nabla {v}^\ve\|_{L^{2}}\|\nabla(\nu_\ve\eta^2)\|_{L^{4}}+\| \zeta_\ve\|_{L^{\infty}}\|\nabla^2(\nu_\ve\eta^2)\|_{L^{\frac43}}+\|\nabla {v}^\ve\|_{L^{2}}+\|\tau_\ve\|_{L^{2}}.
\end{split}
\eeq
Since 
$$
{v}^\ve\eta^2=\omega_\ve\eta^2+\zeta_\ve\nu_\ve\eta^2
$$
we have
\beq
\begin{split}
&\|\nabla^2({v}^\ve\eta^2)\|_{L^{\frac43}}\\
\les&\|\nabla^2(\zeta_\ve\nu_\ve\eta^2)\|_{L^{\frac43}}+\|\nabla^2(\omega_\ve\eta^2)\|_{L^{\frac43}}\\
\les&\| \zeta_\ve\|_{L^{\infty}}\|\nabla^2(\nu_\ve\eta^2)\|_{L^{\frac43}}+\|\nabla \zeta_\ve\|_{L^{2}}\|\nabla (\nu_\ve\eta^2)\|_{L^{4}}+\|\nabla^2\zeta_\ve\nu_\ve\eta^2\|_{L^{\frac43}}+\|\nabla^2(\omega_\ve\eta^2)\|_{L^{\frac43}}\\
\les&\| \zeta_\ve\|_{L^{\infty}}\|\nabla^2(\nu_\ve\eta^2)\|_{L^{\frac43}}+\|\nabla \zeta_\ve\|_{L^{2}}\|\nabla (\nu_\ve\eta^2)\|_{L^{4}}+\|\nabla^2(\zeta_\ve\eta^2)\|_{L^{\frac43}}\\
&+\|\nabla \zeta_\ve\|_{L^{2}}+\|\nabla^2(w_\ve\eta^2)\|_{L^{\frac43}}+1.
\end{split}
\eeq
Therefore, we have 
\beq
\begin{split}
&\|\nabla^2({v}^\ve\eta^2)\|_{L^{\frac43}}\\
\les&\| \zeta_\ve\|_{L^{\infty}}\|\nabla^2(\nu_\ve\eta^2)\|_{L^{\frac43}}+\|\nabla {v}^\ve\|_{L^{2}}\left[1+\|\nabla ({v}^\ve\eta^2)\|_{L^{4}}+\|\nabla (\nu_\ve\eta^2)\|_{L^{4}}\right]+\|\tau_\ve\|_{L^{2}}+1\\
\les&\| \zeta_\ve\|_{L^{\infty}}\|\nabla^2(\nu_\ve\eta^2)\|_{L^{\frac43}}+\|\nabla {v}^\ve\|_{L^{2}}\left[1+\|\nabla ({v}^\ve\eta^2)\|_{L^{4}}\right]+\|\tau_\ve\|_{L^{2}}+1.
\end{split}
\eeq
Since $\nu_\ve=\nu_\ve({v}^\ve)$, we can directly calculate and show that 
\beq
\begin{split}
\|\nabla^2(\nu_\ve\eta^2)\|_{L^{\frac43}}
\les\|\nabla^2({v}^\ve\eta^2)\|_{L^{\frac43}}+\|\nabla {v}^\ve\|_{L^{2}}\left[1+\|\nabla ({v}^\ve\eta^2)\|_{L^{4}}\right]+1.
\end{split}
\eeq
Therefore, we can conclude that
\beq
\begin{split}
\left(1-C\|\zeta_\ve\|_{L^\infty}\right)\|\nabla^2({v}^\ve\eta^2)\|_{L^{\frac43}}
\les\|\nabla {v}^\ve\|_{L^{2}}\left[1+\|\nabla ({v}^\ve\eta^2)\|_{L^{4}}\right]+\|\tau_\ve\|_{L^{2}}+1.
\end{split}
\eeq
By Sobolev's embedding, we have 
\beq
\begin{split}
\|\nabla ({v}^\ve\eta^2)\|_{L^{4}}
\les\|\nabla {v}^\ve\|_{L^{2}}\left[1+\|\nabla ({v}^\ve\eta^2)\|_{L^{4}}\right]+\|\tau_\ve\|_{L^{2}}+1.
\end{split}
\eeq
Taking $\delta_0$ small enough in the assumption \eqref{engdelta}, we conclude that 
\beq
\begin{split}
\|\nabla ({v}^\ve\eta^2)\|_{L^{4}}
\les\|\nabla {v}^\ve\|_{L^{2}}+\|\tau_\ve\|_{L^{2}}+1\leq C(\delta_0, \Lambda_2).
\end{split}
\eeq
Substituting this into (2.19), we obtain 
which implies that 
\beq
\begin{split}
\|\nabla^2 (v_\ve\eta^2)\|_{L^{\frac43}}
\leq C(\delta_0, \Lambda_2).
\end{split}
\eeq
Hence $v^\ve \rightarrow v$ in $H^1(B_{\frac{r_0}3}(\x_0))$. 

By Fubini's theorem, there exists $r_1\in[\frac{r_0}4, \frac{r_0}3]$
\beq\label{estpgs}
\int_{\partial B_{r_1}(\x_0)}|\nabla \zeta_\ve|^2\leq C\int_{B_{\frac{r_0}3}(\x_0)}|\nabla \zeta_\ve|^2\leq C,
\quad \int_{\partial B_{r_1}(\x_0)}|\zeta_\ve|^2\leq C\int_{B_{\frac{r_0}3}(\x_0)}|\zeta_\ve|^2\leq C\ve^2.
\eeq 
Multiplying the equation of $\zeta_\ve$ by $\zeta_\ve$ and integrating by parts over $B_{r_2}$, we obtain
\beq
\int_{B_{r_1}(\x_0)}\big(|\nabla \zeta_\ve|^2+\frac{2}{\ve^2}\chi'(\zeta_\ve)\zeta_\ve^2+|\nabla \nu_\ve|^2\zeta_\ve^2+\nabla \omega_\ve\nabla \nu_\ve \cdot \zeta_\ve\big)-\int_{\partial B_{r_1}(\x_0)}\frac{\partial \zeta_\ve}{\partial\nu} \zeta_\ve=\int_{B_{r_1}(\x_0)}\tau_\ve^{\perp}\zeta_\ve
\eeq
Then we have
\beq
\begin{split}
&\int_{B_{r_1}(\x_0)}\big(|\nabla \zeta_\ve|^2+\frac{2}{\ve^2}\zeta_\ve^2\big)\\
\leq &C\big(\int_{\partial B_{r_1}(\x_0)}|\nabla \zeta_\ve|^2\big)^{\frac12}
\big(\int_{\partial B_{r_1}(\x_0)}|\zeta_\ve|^2\big)^{\frac12}
+C\big(\int_{B_{r_1}(\x_0)}|\nabla \omega_\ve|^{4}\big)^{\frac 12}
\big(\int_{B_{r_1}(\x_0)}|\zeta_\ve|^{4}\big)^{\frac12}\\
&+C\big(\int_{ B_{r_1}(\x_0)}|\tau_\ve|^2\big)^{\frac12}\big(\int_{B_{r_1}(\x_0)}|\zeta_\ve|^2\big)^{\frac12}\leq C\ve.
\end{split}
\eeq
Therefore we have that 
\beq
\frac{\zeta_\ve^2}{\ve^2}\rightarrow 0\ \ \mbox{in }\ \ L^1(B_{r_1}(\x_0)).
\eeq
This completes the proof.
\endpf

\medskip
Now we define the concentration set  by 
\begin{equation}
  \Sigma:=\bigcap_{r>0}\big\{ x\in \Omega:\liminf_{k\to \infty}\int_{B_{r}(\x)}\big(\frac12|\nabla {v}^\ve|^2+\frac{1}{\varepsilon^2}\chi\big(\dist^2({v}^\ve, \mathcal N)\big)\big)> \delta_0^2 \big\},
  \label{}
\end{equation}
where $\delta_0>0$ is  given in Lemma \ref{lemma:biGL}.  We have

\begin{lemma}\label{lemma2.2}  $\Sigma$ is a finite set, and
\beq
v^\ve\rightarrow v\quad\mbox{in }\ \ H^1_{\loc} (\Omega\setminus\Sigma).
\eeq
\end{lemma}
The finiteness of $\Sigma$ follows from a simple covering argument, see also \cite{linwang2016weakthreedim} and \cite{kortum2019concentrationcancellation}.


\section{Convergence of Ginzburg-Landau approximation}

The section is devoted to the proof of Theorem \ref{thm3}. First, recall from the  global energy inequality for \eqref{NGL} that for almost every $t\in(0,T)$, 
\begin{equation}
\begin{split}
  \int_{\Omega\times\left\{ t \right\}}\big( |u^\varepsilon|^2+|\nabla v^\varepsilon|^2+\frac{1}{\varepsilon^2}\chi\big(\dist^2({v}^\ve, \mathcal N)\big)\big)
  +2\int_{Q_t}\big( |\nabla u^\varepsilon|^2+\left| v^\varepsilon_t+u^\varepsilon\cdot \nabla v^\varepsilon\right|^2
  \big)\le E_0.
\end{split}
  \label{GLglobalenergy}
\end{equation}
This, combined with the equation \eqref{NGL}, implies that there exists $p>2$ such that  
\begin{equation}
  \begin{split}
    \sup_{\ve>0}\left[ \left\|u^\ve_t\right\|_{L_t^2 H_x^{-1}+L_t^2 W^{-2, p}}+\left\| v^\ve_t\right\|_{L^{4/3}_t L_x^{4/3}} \right]<\infty.
  \end{split}
  \label{estut}
\end{equation}
Hence, by Aubin-Lions' Lemma,  there exists $u\in L_t^\infty L_x^2 \cap L_t^2 H_x^1(\Omega\times(0,T), \R^2)$ and $v\in L_t^\infty H_x^1\cap L_t^\infty H_x^1(\Omega\times (0,T), \N)$ such that after taking a subsequence,
\begin{equation}\notag
    (u^\ve, v^\ve)\rightarrow (u, v) \text{ in }L^2(\Omega\times(0, T)), \quad
    (\nabla u^\ve, \nabla v^\ve)\rightharpoonup (\nabla u, \nabla v)\text{ in  }L^2(\Omega\times (0,T)).
\end{equation}
Combining this with \eqref{GLglobalenergy}, we obtain 
 $$
 v^\ve_t+u^\ve\cdot \nabla v^\ve\rightharpoonup v_t+u\cdot \nabla v\ \ \text{in }\ \ L^2(\Omega\times (0,T)).
 $$

By the lower semi-continuity, we have 
\begin{equation}
  \begin{split}
    \int_{Q_t}(|\nabla u|^2+|v_t+ u\cdot \nabla v|^2)\le \liminf_{\ve\to 0}\int_{Q_t}(|\nabla u^\ve|^2+|v^\ve_t+u^\ve\cdot \nabla v^\ve|^2)<\infty.
  \end{split}
  \label{}
\end{equation}
By Fatou's Lemma, we have
\begin{equation}
  \begin{split}
    \int_{0}^{t}\liminf_{\ve\to0}\int_\Omega (|\nabla u^\ve|^2+|v^\ve_t+u^\ve\cdot \nabla v^\ve|^2)\le \liminf_{\ve\to 0}\int_{0}^{t}\int_\Omega
    (|\nabla u^\ve|^2+| v^\ve_t+u^\ve\cdot \nabla v^\ve|^2)\leq E_0.
  \end{split}
  \label{}
\end{equation}
Hence there exists $A\subset[0,T]$ with full Lebesgue measure $T$ such that for any $t\in A$ 
\beq
\big(u^\ve(t), v^\ve(t)\big)\rightharpoonup (u(t),v(t))\ \quad\mbox{in }\ \ L^2\times H^1
\eeq
and 
\beq
\liminf\limits_{\ve\rightarrow 0^+}\int_{\Omega}\left(|\nabla u^\ve|^2+|v^\ve_t+u^\ve\cdot\nabla v^\ve|^2\right)(t)<\infty.
\eeq
Now we define the concentration set at $t$ by 
\begin{equation}\label{sigma}
  \Sigma_t:=\bigcap_{r>0}\Big\{ \x\in \Omega:\liminf_{\ve\to 0}\int_{B_r(\x)\times\left\{ t \right\}}\frac{1}{2}|\nabla v^\ve|^2+\frac{1}{\varepsilon^2}\chi\big(\dist^2({v}^\ve, \mathcal N)\big)>\delta_0^2 \Big\},
\end{equation}
where $\delta_0$ is given by Lemma \eqref{lemma:biGL}. 
By Lemma \ref{lemma2.2}, it holds $\#(\Sigma_t)\le C(E_0)$ and
\begin{equation*}
  v^\ve(t)\to v(t) \text{ in }H_{\loc}^1(\Omega\setminus \Sigma(t)).
\end{equation*}

We would first show that $v$ is a weak solution of \eqref{NEL}$_3$ by utilizing the geometric structure as in \cite{ChenStruwe89} (see also \cite{lin2008analysis}). First notice that
there exists a unit vector $\nu_{\N}^\ve\perp T_{\Pi_{\N}(v^\ve)}\N$ such that
\beq\notag
\frac{d}{d{v}}\chi\big(\dist^2({v}^\ve, \N)\big)=2\chi'(\dist^2(v^\ve, \mathcal{N}))\dist({v}^\ve, \N)\nu_{\N}^\ve.
\eeq
Thus for any $\phi\in C^\infty_0(\Omega,\R^L)$ and a.e. $t\in (0,\infty)$ it holds 
\beq\notag
\int_{\Omega\times\{t\}}\langle v^\varepsilon_t+u^\varepsilon\cdot\nabla v^\varepsilon-\Delta v^\varepsilon, 
D\Pi_{\N}(\Pi_{\N}(v^\varepsilon))\phi\rangle=0.
\eeq
If we choose $\phi\in C^\infty_0(\Omega\setminus\Sigma_t)$, 
then it follows from
$\nabla v^\varepsilon\rightarrow \nabla v$ in $H^1_{\rm{loc}}(\Omega\setminus\Sigma_t)$
that, after passing to the limit of the above equation, 
\beq\notag
\int_{\Omega\times\{t\}}\langle v_t+u\cdot\nabla v, D\Pi_{\N}(v)\phi\rangle=
-\int_{\Omega\times\{t\}}\langle\nabla v, \nabla(D\Pi_{\N}(v))\phi\rangle.
\eeq
This implies that 
$$
v_t+u\cdot \nabla v-\Delta v=A_{\N}(v)(\nabla v,\nabla v)
$$
holds weakly in $\Omega\setminus\Sigma_t$. Since $\Sigma_t$ is a finite set, it also holds
weakly in $\Omega$ so that \eqref{NEL}$_3$ holds.

Now, we proceed to verify $u$ satisfies \eqref{NEL}$_1$. First by the estimate \eqref{estut}, we have 
$$u^\ve_t\rightharpoonup u_t,\quad \mbox{in } 
L^2([0,T], H^{-1})\cap L^2([0,T], W^{-2, p})$$
for some $p>2$.
For any $\xi\in C^\infty([0,T])$ with $\xi(T)=0$, $\varphi\in {\bf J}$, since
\beq\notag
\int_{Q_T}u^\ve_t\xi\varphi=-\int_{\Omega}u_0\xi(0)\varphi-\int_{Q_T}u^\ve\xi'\varphi,
\eeq
which, after taking $\epsilon\rightarrow 0$,  implies that 
 $$
 \label{ibpvt}
\int_{Q_T}u_t\xi\varphi=-\int_{\Omega}u_0\xi(0)\varphi-\int_{Q_T}u\xi'\varphi.
$$

\medskip
\noindent{\it Claim:  For any $t\in A$, it holds
\begin{equation}
  \begin{split}
    &0=\int_{\Omega\times\left\{ t \right\}}\left\langle  \pa_t u^\ve, \varphi \right\rangle+\int_{\Omega\times\left\{ t \right\}}\left\langle u^\ve\cdot \nabla u^\ve, \varphi \right\rangle+\int_{\Omega\times \left\{ t \right\}}\left\langle \nabla u^\ve, \nabla \varphi \right\rangle+\int_{\Omega\times\{t\}} (\nabla v^\ve\odot \nabla v^\ve):\nabla \varphi\\
    &\ \ \to \int_{\Omega\times \left\{ t \right\}}\left\langle  u_t, \varphi \right\rangle+\int_{\Omega\times\left\{ t \right\}}\left\langle u\cdot \nabla u, \varphi \right\rangle+\int_{\Omega\times \left\{ t \right\}}\left\langle \nabla u, \nabla \varphi \right\rangle+\int_{\Omega\times\left\{ t \right\}}(\nabla v\odot \nabla v):\nabla \varphi,
  \end{split}
  \label{weakueqn}
\end{equation}
for any $\varphi\in {\bf J}$.
} 
\medskip

 \noindent For this claim, it suffices to show the convergence of  Ericksen stress tensors, i.e., 
\begin{equation}\notag
  \int_{\Omega\times\left\{ t \right\}} (\nabla v^\ve\odot \nabla v^\ve):\nabla \varphi=\int_{\Omega\times\left\{ t \right\}} (\nabla v\odot \nabla v):\nabla \varphi.
  \label{weakcondd}
\end{equation}
For simplicity, we assume $\Sigma_t=\{(0,0)\}\subset\Omega$ consists of a single point at zero. 
Let $\varphi\in C^\infty(\Omega, \R^2)$ be such that 
$\dv \varphi=0$ and $(0,0)\in \mbox{spt}(\varphi)$.  Then
we observe that  by adding $-\frac{1}{2}|\nabla v^\ve|^2\mathbb I_2$, we have
\begin{equation}\notag
  \int_{\Omega\times\left\{ t \right\}} (\nabla v^\ve\odot \nabla v^\ve):\nabla \varphi=\int_{\Omega\times\left\{ t \right\}} \big(\nabla v^\ve\odot \nabla v^\ve-\frac{1}{2}|\nabla v^\ve|^2\mathbb I_2\big):\nabla \varphi.
  \label{weakcondd}
\end{equation}
While by direct computations, we have
\beq\label{2dsplitting}
\nabla v^\ve\odot\nabla v^\ve-\frac{1}{2}|\nabla v^\ve|^2\mathbb I_2
=\frac12\left(
\begin{array}{ll}
|\partial_x v^\ve|^2-|\partial_y v^\ve|^2, &2\langle\partial_x v^\ve,\partial_y v^\ve\rangle\\
2\langle\partial_x v^\ve,\partial_y v^\ve\rangle, &|\partial_y v^\ve|^2-|\partial_x v^\ve|^2
\end{array}
\right).
\eeq
We can assume that there are two real numbers $\alpha, \beta$ such that 
\begin{equation}
  (|\pa_x\vv^\ve|^2-|\pa_y\vv^\ve|^2)dxdy\rightharpoonup (|\pa_x\vv|^2-|\pa_y\vv|^2)dxdy+\alpha \delta_{(0, 0)}, 
\label{convalpha}
\end{equation}
\begin{equation}
  \left\langle \pa_x\vv^\ve, \pa_y\vv^\ve \right\rangle dxdy\rightharpoonup \left\langle \pa_x\vv, \pa_y\vv \right\rangle dxdy+\beta \delta_{(0, 0)},
\label{convbeta}
\end{equation}
hold as convergence of Radon measures. Next we want to show
\begin{equation}\label{vanish}
\alpha=\beta=0.
\end{equation}

Denote
 \begin{equation}
   \Delta\vv^\ve-\frac{1}{\ve^2}\chi'\left( \dist^2(\vv^\ve, \mathcal{N})\right)\frac{d}{d\vv}\left( \dist^2(\vv^\ve, \mathcal{N}) \right)=\f^\ve:= \pa_t\vv^\ve+u^\ve\cdot \nabla\vv^\ve
   \label{}
 \end{equation}
 and
 \begin{equation*}
   e_\ve(\vv^\ve)=\frac{1}{2}|\nabla\vv^\ve|^2+\frac{1}{\ve^2}\chi\left( \dist^2(\vv^\ve, \mathcal{N}) \right).
 \end{equation*}

Now we derive the Pohozaev identity for $v^\epsilon$. For any $X\in C_0^\infty(\Omega,\mathbb R^2)$,
by multiplying the $\vv^\ve$ equation by $X\cdot \nabla \vv^\varepsilon$ and
integrating over  $B_r(0)$ we get
\begin{eqnarray}
&&\int_{\pa B_r(0)}(X^j \vv_j^\varepsilon)\cdot\big( \vv_i^\varepsilon \frac{\x^i}{|\x|} \big)-\int_{B_r(0)}X_i^j \vv_j^\varepsilon \cdot\vv_i^\varepsilon
+\int_{B_r(0)}{\rm{div}} X e_\varepsilon(v^\varepsilon)-\int_{\partial B_r(0)} e_\varepsilon(v^\varepsilon) (X\cdot
\frac{\x}{|\x|})
\nonumber\\
&&=\int_{B_r(0)}(X\cdot \nabla \vv^\varepsilon)\cdot \f^\varepsilon.
\label{poho}
\end{eqnarray}
If we choose $X(\x)=\x$, then we have
\begin{equation}
\begin{split}
r\int_{\pa B_r(0)} \left|\frac{\pa \vv^\varepsilon}{\pa r}\right|^2+\int_{B_r(0)}\frac{2}{\ve^2}\chi(\dist^2(\vv^\ve, \mathcal{N}))
-r\int_{\pa B_r(0)}e_\ve(\vv^\ve)=\int_{B_r(0)}|\x| \frac{\pa \vv^\varepsilon}{\pa r}\cdot \f^\varepsilon .
\end{split}
\label{}
\end{equation}
Then 
\begin{equation}
\begin{split}
\int_{\pa B_r(0)}e_\ve(\vv^\ve)=\int_{\pa B_r(0)}
\left|\frac{\pa \vv^\varepsilon}{\pa r}\right|^2+\frac{1}{r}\int_{B_r(0)}\frac{2}{\ve^2}\chi(\dist^2(\vv^\ve, \mathcal{N})) +
O\big (\int_{B_r(0)}|\nabla \vv^\varepsilon||\f^\varepsilon| \big).
\end{split}
\label{}
\end{equation}
Integrating from $r$ to $R$, we have
\begin{equation}
\begin{split}
\int_{B_R(0)}e_\ve(\vv^\ve)-\int_{B_r(0)}e_\ve(\vv^\ve)&=\int_{B_R(0)\setminus B_r(0)} \left|\frac{\pa \vv^\varepsilon}{\pa r}\right|^2+\int_{r}^{R}\frac{1}{\tau}\int_{B_\tau(0)}\frac{2}{\ve^2}\chi(\dist^2(\vv^\ve, \mathcal{N}))d\tau\\
&+\int_{r}^{R}{O}\big( \int_{B_\tau(0)}|\nabla \vv^\varepsilon||\f^\varepsilon| \big)d\tau.
\end{split}
\label{5.10}
\end{equation}
Since $\Sigma_t=\{(0,0)\}$, we can assume that 
\begin{equation}
e_\ve(\vv^\ve)d\x\rightharpoonup \frac12|\nabla\vv|^2 d\x+\gamma \delta_{(0, 0)}, \qquad \text{ in }\ \ B_\delta(0)
\label{}
\end{equation}
as convergence of Radon measures,  where $\gamma\ge 0$. Since $t\in A$, 
\begin{align*}
\lim_{\varepsilon\to 0}\int_{B_\tau(0)}|\f^\varepsilon||\nabla \vv^\varepsilon|&\le \lim_{\varepsilon\to 0}\big( \int_{B_\tau(0)}|\f^\varepsilon|^2 \big)^{\frac{1}{2}}\big( \int_{B_\tau(0)}|\nabla \vv^\varepsilon|^2 \big)^\frac{1}{2}\le C E_0,
\end{align*}
Hence, by sending $\varepsilon\to 0$ we obtain from \eqref{5.10} that
\begin{align*}
\int_{B_R(0)\setminus B_r(0)} \frac12|\nabla\vv|^2&\ge \int_{B_R(0)\setminus B_r(0)}\left|\frac{\pa\vv}{\pa r}\right|^2
+\int_{r}^{R}\frac{1}{\tau}\lim_{\ve\to0}\int_{B_\tau(0)}\frac{2}{\ve^2}\chi(\dist^2(\vv^\ve, \mathcal{N})) d\tau
+{O}(R).
\end{align*}
Sending $r\to 0$, we have
\begin{align*}
\int_{B_R(0)} \frac{1}{2}|\nabla\vv|^2&\ge \int_{B_R(0)} \left|\frac{\pa \vv}{\pa r}\right|^2
+\int_{0}^{R}\frac{1}{\tau}
\lim_{\ve\to0}\int_{B_\tau(0)}\frac{2}{\ve^2}\chi(\dist^2(\vv^\ve, \mathcal{N})) d\tau+{O}(R).
\end{align*}
From this, we claim that 
\begin{equation}
\frac{2}{\ve^2}\chi(\dist^2(\vv^\ve, \mathcal{N}))\to 0\ \  \text{ in }\ \ \ L^1(B_\delta).
\label{biaxialL1strongconv}
\end{equation}
For, otherwise, 
\begin{equation*}
\frac{2}{\ve^2}\chi(\dist^2(\vv^\ve, \mathcal{N}))\,dx\rightharpoonup  \kappa \delta_{(0,0)}
\end{equation*}
for some $\kappa>0$, this implies
\begin{equation*}
\int_{0}^{R}\frac{1}{\tau}\lim_{\ve\to0}\int_{B_\tau}\frac{2}{\ve^2}\chi(\dist^2(\vv^\ve, \mathcal{N}))=\int_{0}^{R}\frac{\kappa}{\tau}d\tau=\infty,
\end{equation*}
which is impossible. 

Choosing $X(\x)=(x, 0)$ in \eqref{poho}, we obtain that
\begin{eqnarray}\label{hopf1}
&&\frac12\int_{B_r(0)}\Big( \big|{\pa_y \vv^\ve}\big|^2-\big|{\pa_x \vv^\ve}\big|^2 \Big)
+\int_{B_r(0)} \frac{1}{\varepsilon^2}\chi(\dist^2(\vv^\ve, \mathcal{N}))\nonumber\\
&&=\int_{B_r(0)} x\langle{\partial_x v^\ve}, \f^\ve\rangle +\int_{\partial B_r(0)} \frac{x^2}{r} e_\ve(v^\ve)
-\int_{\partial B_r(0)} x\langle{\partial_x v^\ve}, \frac{\pa v^\ve}{\pa r}\rangle.
\end{eqnarray}
Observe that by Fubini's theorem, for a.e. $r>0$ it holds that 
\begin{equation*}
\begin{split}
\int_{\partial B_r(0)} x\langle{\partial_x v^\ve}, \frac{\pa v^\ve}{\pa r}\rangle
&\to \int_{\pa B_r(0)} x\langle{\pa_x \vv}, \frac{\pa \vv}{\pa r}\rangle, \\
\int_{\pa B_r(0)}\frac{x^2}{r}e_\varepsilon(\vv^\ve) &\to \frac{1}{2}\int_{\pa B_r}\frac{x^2}{r}|\nabla \vv|^2,
\end{split}
\end{equation*}
and by \eqref{biaxialL1strongconv}, 
\begin{equation*}
\int_{B_r(0)}\frac{1}{\varepsilon^2}\chi(\dist^2(\vv^\ve, \mathcal{N}))\to 0.
\end{equation*}
Furthermore,
\begin{align*}
\big|\int_{B_r(0)}x\langle{\pa_x \vv^\ve},  \f^\ve\rangle\big|
&\le Cr\left\|\f^\ve\right\|_{L^2}\left\|\nabla \vv^\ve\right\|_{L^2}={O}(r).
\end{align*}
Hence, by sending $\varepsilon\to 0$ in \eqref{hopf1}, we obtain 
\begin{equation*}
\int_{B_r(0)} \big( \big|{\pa_y \vv}\big|^2-\big|{\pa_x \vv}\big|^2 \big)+\alpha={O}(r),
\end{equation*}
this further implies $\alpha=0$ after sending $r\to 0$.

Similarly, if we choose $X(\x)=(0, x)$ in \eqref{poho} and pass the limit in
 the resulting equation, we can  get that 
\begin{equation*}
\int_{B_r(0)}\big\langle {\pa_x \vv}, {\pa_y \vv}\big\rangle+\beta ={O}(r).
\end{equation*}
Hence $\beta=0$. This proves \eqref{vanish} and hence completes the proof of Claim.

\medskip
Multiplying \eqref{weakueqn} by $\xi\in C^\infty([0,T])$ with $\xi(T)=0$ and integrating over $[0,T]$, we conclude that $u$ satisfies the \eqref{NEL}$_1$ on $Q_T$.
The proof of Theorem \ref{thm3} is complete.


\section{Compactness of simplified Ericksen-Leslie system}

This section is devoted to prove Theorem \ref{thm1}. First notice that since the sequence of weak solutions $(u^k, v^k)$
satisfies the assumption \eqref{kengbd}, and 
\beq\notag
(u^k_0,v^k_0)\rightharpoonup (u_0, v_0),\quad \mbox{in } L^2(\Omega)\times H^1(\Omega), 
\eeq
 there exists $(u(\x,t), v(\x,t)):\Omega\times (0,T)\rightarrow \R^2\times\N$ such that 
\beq\label{engcon1}
(u^k, v^k)\rightharpoonup (u,v) \ \mbox{in }\ L^2([0,T],H^1(\Omega)),
\eeq
\beq\label{engcon2}
v^k_t+u^k\cdot \nabla v^k\rightharpoonup v_t+u\cdot \nabla v\ \ \text{in }\ L^2([0,T], L^2(\Omega)).
\eeq
Also it follows from  \eqref{NEL} and \eqref{kengbd} that there exists $p>2$ that 
\begin{equation}
  \begin{split}
    &\sup_{k} \Big[\big\|u^k_t\big\|_{L_t^2 H_x^{-1}+L_t^2 W^{-2, p}_x}+\big\|v^k_t\big\|_{L_t^{2}H_x^{-1}}\Big]<\infty.
  \end{split}
  \label{}
\end{equation}
Hence, by Aubin-Lions' Lemma we have that 
\begin{equation*}
  (u^k, v^k)\to (u, v) \text{ in } L^2(Q_T)\times L^2(Q_T). 
\end{equation*}

By the lower semi-continuity, we have
\beq\notag
\begin{split}
\int_{Q_t}\left(|\nabla u|^2+|v_t+u\cdot\nabla v|^2\right)
\leq \liminf\limits_{k\rightarrow\infty}\int_{Q_t}\left(|\nabla u^k|^2+|v^k_t+u^k\cdot\nabla v^k|^2\right)\leq C_0.
\end{split}
\eeq
By Fatou's Lemma and \eqref{kengbd}, we have 
\beq\notag
\begin{split}
\int_0^t\liminf\limits_{k\rightarrow\infty}\int_{\Omega}\left(|\nabla u^k|^2+|v^k_t+u^k\cdot\nabla v^k|^2\right)
\leq \liminf\limits_{k\rightarrow\infty}\int_{Q_t}\left(|\nabla u^k|^2+|v^k_t+u^k\cdot\nabla v^k|^2\right)\leq C_0.
\end{split}
\eeq
Hence, there exists $A\subset[0,T]$ with full Lebesgue measure $T$, such that for all $t\in A$ 
\beq
\big(u^k(t), v^k(t)\big)\rightharpoonup (u(t),v(t)),\quad\mbox{in }L^2(\Omega)\times H^1(\Omega)
\eeq
and 
\beq
\liminf\limits_{k\rightarrow\infty}\int_{\Omega}\left(|\nabla u^k|^2+|v^k_t+u^k\cdot\nabla v^k|^2\right)(t)<\infty.
\eeq
Now we define the concentration set  at time $t\in (0,T]$ by 
\begin{equation}
  \Sigma_t:=\bigcap_{r>0}\Big\{ x\in \Omega:\liminf_{k\to \infty}\int_{B_r(x)}|\nabla v^k|^2> \delta_0^2 \Big\}
  \label{}
\end{equation}
where $\delta_0$ is small constant given by Theorem 1.2 in \cite{SharpTopping2013}. 
As in \cite{SharpTopping2013} (see also \cite{linwang2016weakthreedim}, \cite{wang1999remark}), 
we can show that for any $t\in A$, it holds that $\#(\Sigma_t)\leq C(E_0)$ and 
\beq
v^k(t)\rightarrow v\ \quad\mbox{in }\ \ H^1_{\loc} (\Omega\setminus\Sigma_t).
\eeq

Similar to the proof of Theorem \ref{thm3}, we can show the weak limit $(u,v)$ satisfies the third equation of \eqref{NEL} in the weak sense. It remains to show that the first equation of \eqref{NEL} is also valid in the weak sense. 

Similar to the proof of Theorem \ref{thm3}, to complete the proof of Theorem \ref{thm1}, it is suffices to show
\beq\label{H1conv}
\lim\limits_{k\rightarrow\infty}\int_{\Omega\times\{t\}}\left(\nabla v^k\odot\nabla v^k\right):\nabla \varphi=\int_{\Omega\times\{t\}}\left(\nabla v\odot\nabla v\right):\nabla \varphi, \ \forall\varphi\in \mathbf{J}.
\eeq
For simplicity, assume $\Sigma_t=\{(0,0)\}\subset\Omega$.
Let $\varphi\in C^\infty(\Omega, \R^2)$ be such that $\dv \varphi=0$ and $(0,0)\in \mbox{spt}(\varphi)$. 
By the same calculation as in \eqref{2dsplitting}, we have
\beq\notag
\nabla v^\ve\odot\nabla v^\ve-\frac{1}{2}|\nabla v^\ve|^2\mathbb I_2
=\frac12\left(
\begin{array}{ll}
|\partial_x v^\ve|^2-|\partial_y v^\ve|^2, &2\langle\partial_x v^\ve,\partial_y v^\ve\rangle\\
2\langle\partial_x v^\ve,\partial_y v^\ve\rangle, &|\partial_y v^\ve|^2-|\partial_x v^\ve|^2
\end{array}
\right).
\eeq

For any $t\in A$, $v^k(t)$ is an approximated harmonic maps from $\Omega$ to $\N$:
\begin{equation}
\Delta v^k(t)+A(v^k)(\nabla v^k,\nabla v^k)=g^k(t):=v^k_t(t)+u^k\cdot\nabla v^k(t)\in L^2(\Omega).
\end{equation}
Recall the Hopf differential of $v^k$ is defined by
\beq
\mathcal H^k=\big(\frac{\partial v^k}{\partial z}\big)^2
=\big|{\partial _xv^k}\big|^2-\big|{\partial_y v^k}\big|^2
+2i\big\langle{\partial_x v^k}, {\partial_y v^k}\big\rangle,
\eeq
where $z=x+iy\in \mathbb C$.
Then 
\beq
\frac{\partial {{\mathcal{H}}}^k}{\partial\bar z}=2\frac{\partial_x v^k}{\partial z}\frac{\partial^2 v^k}{\partial \bar z\partial z}
=2\Delta v^k\frac{\partial v^k}{\partial z}
=2g^k(t)\cdot\frac{\partial v^k}{\partial z}:=G^k.
\eeq
It is clear that 
\beq
\|G^k\|_{L^1(B_r)}\leq 2\|g^k(t)\|_{L^2(\Omega)}\big\|\frac{\partial v^k}{\partial z}\big\|_{L^2(\Omega)}
\leq 2C_0.
\eeq
Therefore, for any $z\in B_r(0)$
\beq
{{\mathcal{H}}}^k(z)=\int_{\partial B_{2r}(0)}\frac{{{\mathcal{H}}}^k(\omega)}{z-\omega}\,d\sigma+\int_{B_{2r}(0)}\frac{G^k(\omega)}{z-\omega}\,d\omega.
\eeq
By the Young inequality of convolutions, we obtain
\beq
\|{{\mathcal{H}}}^k\|_{L^p(B_r)}
\leq C(r, p) \|\mathcal{H}^k\|_{L^1(\partial B_{2r})}+\big\|\frac{1}{z}\big\|_{L^p}\|G^k\|_{L^1( B_{2r})}\leq C(r, p).
\eeq
for any $1<p<2$. 
From this, we immediately conclude that
\beq\notag
|\partial_x v^k|^2-|\partial_yv^k|^2\rightharpoonup
|\partial_x v|^2-|\partial_yv|^2, \quad \langle\partial_x v^k,\partial_yv^k\rangle\rightharpoonup  \langle\partial_x v,\partial_y v\rangle\ \quad \mbox{in }\ \ L^p(B_r(0))
\eeq
for any $1<p<2$, which implies \eqref{H1conv}. 
This completes the proof of Theorem \ref{thm1}.

\bigskip
\bigskip
\noindent{\bf Acknowledgments.} The first and third author are partially supported by NSF DMS 1764417.



\begin{thebibliography}{99}
\bibitem{ChenStruwe89}
Yun~Mei Chen and Michael Struwe, \emph{Existence and partial regularity results
	for the heat flow for harmonic maps}, Math. Z. \textbf{201} (1989), no.~1,
83--103.

\bibitem{DipernaMajda1988concentration}
Ronald~J. DiPerna and Andrew Majda, \emph{Reduced {H}ausdorff dimension and
	concentration-cancellation for two-dimensional incompressible flow}, J. Amer.
Math. Soc. \textbf{1} (1988), no.~1, 59--95.

\bibitem{GoversVertogen1984a} E. Govers, G. Vertogen, {\it Elastic continuum theory of biaxial nematics}.
Phys. Rev. A, Vol. 30, No. 4 (1984).

\bibitem{GoversVertogen1984b} E. Govers, G. Vertogen, {\it Erratum: Elastic continuum theory of biaxial nematics}
[Phys. Rev. A, Vol. 30, No. 4 (1984)], Phys. Rev. A, Vol. 31, No. 3 (1985).

\bibitem{GoversVertogen1985} E. Govers, G. Vertogen,  {\it FLUID DYNAMICS OF BIAXIAL NEMATICS}.
Physica 133A (1985) 337-344.


\bibitem{Huang2016FiniteTS}
Tao Huang, Fanghua Lin, Chun Liu, and Changyou Wang, \emph{Finite time
	singularity of the nematic liquid crystal flow in dimension three}, Arch. Ration. Mech. Anal. \textbf{221} (2016), 1223--1254.

\bibitem{kortum2019concentrationcancellation}
Joshua Kortum, \emph{Concentration-cancellation in the {E}ricksen-{L}eslie
	model}, arXiv:1912.10429 (2019).

\bibitem{Krylov2008Sobolev}
N.~V. Krylov, \emph{Lectures on elliptic and parabolic equations in {S}obolev
	spaces}, Graduate Studies in Mathematics, vol.~96, American Mathematical
Society, Providence, RI, 2008.

\bibitem{Hong2011} MinChun  Hong, {\em Global existence of solutions of the simplified Ericksen-Leslie system in dimension two.} Calc. Var. Partial Differential Equations 40, no. 1-2 (2011), 15-36.

\bibitem{LLWWZ19}
Chen-Chih Lai, Fanghua Lin, Changyou Wang, Juncheng Wei, and Yifu Zhou,
\emph{Finite time blow-up for the nematic liquid crystal flow in dimension
	two}, arXiv:1908.10955 (2019).

\bibitem{lin1989nonlinear}
Fanghua Lin, \emph{Nonlinear theory of defects in nematic liquid crystals;
	phase transition and flow phenomena}, Comm. Pure Appl. Math. \textbf{42} (1989), no.~6, 789--814.

\bibitem{lin2010liquid}
Fanghua Lin, Junyu Lin, and Changyou Wang, \emph{Liquid crystal flows in two
	dimensions}, Arch. Ration. Mech. Anal. \textbf{197} (2010),
no.~1, 297--336.

\bibitem{linliu1995nonparabolic}
Fanghua Lin and Chun Liu, \emph{Nonparabolic dissipative systems modeling the
	flow of liquid crystals}, Comm. Pure Appl. Math. \textbf{48} (1995), no.~5,
501--537.

\bibitem{linliu1996partialregularity}
Fanghua Lin and Chun Liu, \emph{Partial regularity of the dynamic system modeling the flow of
	liquid crystals}, Discrete Contin. Dynam. Systems \textbf{2} (1996), no.~1,
1--22.

\bibitem{lin2008analysis}
Fanghua Lin and Changyou Wang, \emph{The analysis of harmonic maps and their
	heat flows}, World Scientific, 2008.

\bibitem{lin2010uniqueness}
Fanghua Lin and Changyou Wang, \emph{On the uniqueness of heat flow of harmonic maps and hydrodynamic
	flow of nematic liquid crystals}, Chinese Annals of Mathematics, Series B
\textbf{31} (2010), no.~6, 921--938.

\bibitem{lin2014recent}
Fanghua Lin and Changyou Wang, \emph{Recent developments of analysis for hydrodynamic flow of nematic
	liquid crystals}, Philosophical Transactions of the Royal Society A:
Mathematical, Physical and Engineering Sciences \textbf{372} (2014),
no.~2029, 20130361.

\bibitem{linwang2016weakthreedim}
Fanghua Lin and Changyou Wang, \emph{Global existence of weak solutions of the nematic liquid crystal
	flow in dimension three}, Comm. Pure Appl. Math. \textbf{69} (2016), no.~8,
1532--1571.

\bibitem{LinLiWang20}
Junyu Lin, Yimei Li, and Changyou Wang, \emph{On static and hydrodynamic
	biaxial nematic liquid crystals}. Preprint (2020).

\bibitem{MajdaBertozzi2002Vorticity}
Andrew~J. Majda and Andrea~L. Bertozzi, \emph{Vorticity and incompressible
	flow}, Cambridge Texts in Applied Mathematics, vol.~27, Cambridge University
Press, Cambridge, 2002.

\bibitem{SharpTopping2013}
Ben Sharp and Peter Topping, \emph{Decay estimates for {R}ivi\`ere's equation,
	with applications to regularity and compactness}, Trans. Amer. Math. Soc.
\textbf{365} (2013), no.~5, 2317--2339.

\bibitem{wang1999remark}
Changyou Wang, \emph{A remark on harmonic map flows from surfaces},
Differential Integral Equations \textbf{12} (1999), no.~2, 161--166.	

	
\end{thebibliography}
\end{document}